\input vanilla.sty
\scaletype{\magstep1}
\scalelinespacing{\magstep1}
\def\bull{\vrule height .9ex width .8ex depth -.1ex}

\title On vector-valued inequalities for Sidon sets and sets
of interpolation \endtitle

\author N.J.  Kalton \footnote{Research supported by
NSF-grant DMS-8901636}\\ University of Missouri-Columbia
\endauthor

\vskip1truecm

\subheading{Abstract} Let $E$ be a Sidon subset of the
integers and suppose $X$ is a Banach space.  Then Pisier has
shown that $E$-spectral polynomials with values in $X$
behave like Rademacher sums with respect to $L_p-$norms.  We
consider the situation when $X$ is a quasi-Banach space.
For general quasi-Banach spaces we show that a similar
result holds if and only if $E$ is a set of interpolation
($I_0$-set).  However for certain special classes of
quasi-Banach spaces we are able to prove such a result for
larger sets.  Thus if $X$ is restricted to be ``natural''
then the result holds for all Sidon sets.  We also consider
spaces with plurisubharmonic norms and introduce the class
of analytic Sidon sets.

\vskip2truecm

\subheading{1.  Introduction}

Suppose $G$ is a compact abelian group.  We denote by
$\mu_G$ normalized Haar measure on $G$ and by $\Gamma$ the
dual group of $G$.  We recall that a subset $E$ of $\Gamma$
is called a {\it Sidon set} if there is a constant $M$ such
that for every finitely nonzero map $a:E\to\bold C$ we have
$$ \sum_{\gamma\in E}|a(\gamma)|\le M\max_{g\in
G}|\sum_{\gamma\in E}a(\gamma)\gamma(g)|.$$

We define $\Delta$ to be the Cantor group i.e.
$\Delta=\{\pm 1\}^{\bold N}$.  If $t\in \Delta$ we denote by
$\epsilon_n(t)$ the $n$-th co-ordinate of $t.$ The sequence
$(\epsilon_n)$ is an example of a Sidon set.  Of course the
sequence $(\epsilon_n)$ is a model for the Rademacher
functions on $[0,1].$ Similarly we denote the co-ordinate
maps on $\bold T^{\bold N}$ by $\eta_n.$

Suppose now that $G$ is a compact abelian group.  If $X$ is
a Banach space, or more generally a quasi-Banach space with
a continuous quasinorm and $\phi:G\to X$ is a Borel map we
define $\|\phi\|_p$ for $0<p\le\infty$ to be the $L_p-$norm
of $\phi$ i.e.  $\|\phi\|_p=(\int_G
\|\phi(g)\|^pd\mu_G(g))^{1/p}$ if $0<p<\infty$ and
$\|\phi\|_{\infty}=\text{ess }\sup_{g\in G}\|\phi(g)\|.$

It is a theorem of Pisier [12] that if $E$ is a Sidon set
then there is a constant $M$ so that for every subset
$\{\gamma_1,\ldots,\gamma_n\}$ of $E$, every
$x_1,\ldots,x_n$ chosen from a Banach space $X$ and every
$1\le p\le\infty$ we have $$
M^{-1}\|\sum_{k=1}^nx_k\epsilon_k\|_p \le \|\sum_{k=1}^n
x_k\gamma_k\|_p \le M\|\sum_{k=1}^nx_k\epsilon_k\|_p.  \tag
*$$ Thus a Sidon set behaves like the Rademacher sequence
for Banach space valued functions.  The result can be
similarly stated for $(\eta_n)$ in place of $(\epsilon_n).$
Recently Asmar and Montgomery-Smith [1] have taken Pisier's
ideas further by establishing distributional inequalities in
the same spirit.

It is natural to ask whether Pisier's inequalities can be
extended to arbitrary quasi-Banach spaces.  This question
was suggested to the author by Asmar and Montgomery-Smith.
For convenience we suppose that every quasi-Banach space is
$r$-normed for some $r<1$ i.e. the quasinorm satisfies
$\|x+y\|^r\le \|x\|^r+\|y\|^r$ for all $x,y$; an $r$-norm is
necessarily continuous.  We can then ask, for fixed $0<p\le
\infty$ for which sets $E$ inequality (*) holds, if we
restrict $X$ to belong to some class of quasi-Banach spaces,
for some constant $M=M(E,X).$

It turns out Pisier's results do not in general extend to
the non-locally convex case.  In fact we show that if we fix
$r<1$ and ask that a set $E$ satisfies (*) for some fixed
$p$ and every $r$-normable quasi-Banach space $X$ then this
condition precisely characterizes sets of interpolation as
studied in [2], [3], [4], [5], [8], [9], [13] and [14].  We
recall that $E$ is called a {\it set of interpolation (set
of type ($I_0$))} if it has the property that every
$f\in\ell_{\infty}(E)$ (the collection of all bounded
complex functions on $E$) can be extended to a continuous
function on the Bohr compactification $b\Gamma$ of $\Gamma.$

However, in spite of this result, there are specific classes
of quasi-Banach spaces for which (*) holds for a larger
class of sets $E.$ If we restrict $X$ to be a natural
quasi-Banach space then (*) holds for all Sidon sets $E.$
Here a quasi-Banach space is called {\it natural} if it is
linearly isomorphic to a closed linear subspace of a
(complex) quasi-Banach lattice $Y$ which is $q$-convex for
some $q>0,$ i.e. such that for a suitable constant $C$ we
have $$ \|(\sum_{k=1}^n |y_k|^q)^{1/q}\| \le C(\sum_{k=1}^n
\|y_k\|^q)^{1/q}$$ for every $y_1,\ldots,y_n\in Y.$ Natural
quasi-Banach spaces form a fairly broad class including
almost all function spaces which arise in analysis.  The
reader is referred to [6] for a discussion of examples.
Notice that, of course, the spaces $L_q$ for $q<1$ are
natural so that, in particular, (*) holds for all $p$ and
all Sidon sets $E$ for every $0<p\le\infty.$ The case $p=q$
here would be a direct consequence of Fubini's theorem, but
the other cases, including $p=\infty$ are less obvious.

A quasi-Banach lattice $X$ is natural if and only if it is
{\it A-convex}, i.e. it has an equivalent plurisubharmonic
quasi-norm.  Here a quasinorm is plurisubharmonic if it
satisfies $$ \|x\| \le \int_0^{2\pi}
\|x+e^{i\theta}y\|\frac{d\theta}{2\pi}$$ for every $x,y\in
X.$ There are examples of A-convex spaces which are not
natural, namely the Schatten ideals $S_p$ for $p<1$ [7].  Of
course, it follows that $S_p$ cannot be embedded in any
quasi-Banach lattice which is A-convex when $0<p<1.$ Thus we
may ask for what sets $E$ (*) holds for every A-convex
space.  Here, we are unable to give a precise
characterization of the sets $E$ such that (*) holds.  In
fact we define $E$ to be an analytic Sidon set if (*) holds,
for $p=\infty$ (or, equivalently for any other
$0<p<\infty)$, for every A-convex quasi-Banach space $X$.
We show that any finite union of Hadamard sequences in
$\bold N\subset \bold Z$ is an analytic Sidon set.  In
particular a set such as $\{3^n\}\cup\{3^n+n\}$ is an
analytic Sidon set but not a set of interpolation.  However
we have no example of a Sidon set which is not an analytic
Sidon set.

We would like to thank Nakhle Asmar, Stephen
Montgomery-Smith and David Grow for their helpful comments
on the content of this paper.

\vskip1truecm

\subheading{2.  The results}

Suppose $G$ is a compact abelian group and $\Gamma$ is its
dual group.  Let $E$ be a subset of $\Gamma.$ Suppose $X$ is
a quasi-Banach space and that $0<p\le\infty$; then we will
say that $E$ has property $\Cal C_p(X)$ if there is a
constant $M$ such that for any finite subset
$\{\gamma_1,\ldots,\gamma_n\}$ of $E$ and any
$x_1,\ldots,x_n$ of $X$ we have $(*)$ i.e.  $$
M^{-1}\|\sum_{k=1}^nx_k\epsilon_k\|_p\le
\|\sum_{k=1}^nx_k\gamma_k\|_p\le
M\|\sum_{k=1}^nx_k\epsilon_k\|_p.$$ (Note that in contrast
to Pisier's result (*), we here assume $p$ fixed.)  We start
by observing that $E$ is a Sidon set if and only if $E$ has
property $\Cal C_{\infty}(\bold C).$ It follows from the
results of Pisier [12] a Sidon set has property $\Cal
C_p(X)$ for every Banach space $X$ and for every
$0<p<\infty.$ See also Asmar and Montgomery-Smith [1] and
Pelczynski [11].

Note that for any $t\in \Delta$ we have that
$\|\sum\epsilon_k(t)x_k\epsilon_k\|_p=\|\sum
x_k\epsilon_k\|_p.$ Now any real sequence $(a_1,\ldots,a_n)$
with $\max |a_k|\le 1$ can be written in the form
$a_k=\sum_{j=1}^{\infty}2^{-j}\epsilon_k(t_j)$ and it
follows quickly by taking real and imaginary parts that the
there is a constant $C=C(r,p)$ so that for any complex
$a_1,\ldots,a_n$ and any $r$-normed space $X$ we have $$
\|\sum_{k=1}^n a_kx_k\epsilon_k\|_p\le
C\|a\|_{\infty}\|\sum_{k=1}^n x_k\epsilon_k\|_p.$$ From this
it follows quickly that $\|\sum_{k=1}^nx_k\eta_k\|_p$ is
equivalent to $\|\sum_{k=1}^nx_k\epsilon_k\|_p.$ In
particular we can replace $\epsilon_k$ by $\eta_k$ in the
definition of property $\Cal C_p(X)$.

We note that if $E$ has property $\Cal C_p(X)$ then it is
immediate that $E$ has property $\Cal C_p(\ell_p(X))$ and
further that $E$ has property $\Cal C_p(Y)$ for any
quasi-Banach space finitely representable in $X$ (or, of
course, in $\ell_p(X))$.

For a fixed quasi-Banach space $X$ and a fixed subset $E$ of
$\Gamma$ we let $\Cal P_E(X)$ denote the space of $X$-valued
$E$-polynomials i.e. functions $\phi:G\to X$ of the form
$\phi=\sum_{\gamma\in E} x(\gamma)\gamma$ where $x(\gamma)$
is only finitely nonzero.  If $f\in \ell_{\infty}(E)$ we
define $T_f:\Cal P_E(X)\to\Cal P_E(X)$ by $$ T_f(\sum
x(\gamma)\gamma) =\sum f(\gamma)x(\gamma)\gamma.$$ We then
define $\|f\|_{\Cal M_p(E,X)}$ to the operator norm of $T_f$
on $\Cal P_E(X)$ for the $L_p-$norm (and to be $\infty$ is
this operator is unbounded).

\proclaim{Lemma 1}In order that $E$ has property $\Cal
C_p(X)$ it is necessary and sufficient that there exists a
constant $C$ such that $$ \|f\|_{\Cal M_p(E,X)}\le
C\|f\|_{\infty}$$ for all
$f\in\ell_{\infty}(E).$\endproclaim

\demo{Proof}If $E$ has property $\Cal C_p(X)$ then it also
satisfies (*) for $(\eta_n)$ in place of $(\epsilon_n)$ for
a suitable constant $M.$ Thus if $f\in\ell_{\infty}(E)$ and
$\phi\in \Cal P_E(X)$ then $$ \|T_f\phi\|_p \le
M^2\|f\|_{\infty}\|\phi\|_p.$$

For the converse direction, we consider the case $p<\infty.$
Suppose $\{\gamma_1,\ldots,\gamma_n\}$ is a finite subset of
$E.$ Then for any $x_1,\ldots,x_n$ $$ \align
C^{-p}\int_{\bold T^{\bold
N}}\|\sum_{k=1}^nx_k\eta_k\|^pd\mu_{\bold T^{\bold N}} &=
C^{-p}\int_{\bold T^{\bold N}}\int_G
\|\sum_{k=1}^nx_k\eta_k(s)\gamma_k(t)\|^pd\mu_{\bold
T^{\bold N}}(s) d\mu_G(t)\\ &\le
\int_G\|\sum_{k=1}^nx_k\gamma_k\|^pd\mu_G\\ &\le
C^{p}\int_{\bold T^{\bold N}}\int_G
\|\sum_{k=1}^nx_k\eta_k(s)\gamma_k(t)\|^pd\mu_{\bold
T^{\bold N}}(s) d\mu_G(t)\\ &\le C^p\int_{\bold T^{\bold
N}}\|\sum_{k=1}^nx_k\eta_k\|^pd\mu_{\bold T^{\bold N}}.
\endalign $$ This estimate together with a similar estimate
in the opposite direction gives the conclusion.  The case
$p=\infty$ is similar.  \bull\enddemo

If $E$ is a subset of $\Gamma$, $N\in\bold N$ and $\delta>0$
we let $AP(E,N,\delta)$ be the set of $f\in\ell_{\infty}(E)$
such that there exist $g_1,\ldots,g_N\in G$ (not necessarily
distinct) and $\alpha_1,\ldots,\alpha_N\in\bold C$ with
$\max_{1\le j\le N}|\alpha_j|\le 1$ and $$ |f(\gamma) -
\sum_{j=1}^N \alpha_j\gamma(g_j)|\le \delta$$ for $\gamma\in
E.$

The following theorem improves slightly on results of Kahane
[5] and M\'ela [8].  Perhaps also, our approach is slightly
more direct.  We write
$B_{\ell_{\infty}(E)}=\{f\in\ell_{\infty}(E):\
\|f\|_{\infty}\le 1\}.$

\proclaim{Theorem 2}Let $G$ be a compact abelian group and
let $\Gamma$ be its dual group.  Suppose $E$ is a subset of
$\Gamma$.  Then the following conditions on $E$ are
equivalent:\newline (1) $E$ is a set of
interpolation.\newline (2) There exists an integer $N$ so
that $B_{\ell_{\infty}(E)}\subset AP(E,N,1/2).$\newline (3)
There exists $M$ and $0<\delta<1$ so that if $f\in
B_{\ell_{\infty}(E)}$ then there exist complex numbers
$(c_j)_{j=1}^{\infty}$ with $|c_j|\le M\delta^j$ and
$(g_j)_{j=1}^{\infty}$ in $G$ with $$ f(\gamma) =
\sum_{j=1}^{\infty}c_j\gamma(g_j)$$ for $\gamma\in
E.$\endproclaim

\demo{Proof}$(1)\Rightarrow (2).$ It follows from the
Stone-Weierstrass Theorem that $$\bold T^E\subset
\cup_{m=1}^{\infty}AP(E,m,1/5).$$ Let $\mu=\mu_{\bold T^E}$.
Since each $AP(E,m,1/5)\cap \bold T^E$ is closed it is clear
that there exists $m$ so that $\mu(AP(E,m,1/5)\cap \bold
T^E)>1/2.$ Thus if $f\in \bold T^E$ we can find $f_1,f_2\in
AP(E,m,1/5)\cap \bold T^E$ so that $f=f_1f_2.$ Hence $f\in
AP(E,m^2,1/2).$ This clearly implies (2) with $N=2m^2.$

$(2)\Rightarrow (3).$ We let $\delta=2^{-1/N}$ and $M=2.$
Then given $f\in B_{\ell_{\infty}(E)}$ we can find
$(c_j)_{j=1}^N$ and $(g_j)_{j=1}^N$ with $|c_j|\le 1\le
M\delta^{j}$ and $$ |f(\gamma) - \sum_{j=1}^N c_j
\gamma(g_j)| \le 1/2 $$ for $\gamma\in E.$ Let
$f_1(\gamma)=2(f(\gamma)-\sum_{j=1}^Nc_j\gamma(g_j))$ and
iterate the argument.

$(3)\Rightarrow (1)$:  Obvious.\bull\enddemo

\proclaim{Theorem 3} Suppose $G$ is a compact abelian group,
$E$ is a subset of the dual group $\Gamma$ and that $0<r<1,\
0<p\le\infty.$ In order that $E$ satisfies $\Cal C_p(X)$ for
every $r$-normable quasi-Banach space $X$ it is necessary
and sufficient that $E$ be a set of
interpolation.\endproclaim

\demo{Proof}First suppose that $E$ is a set of interpolation
so that it verifies (3) of Theorem 2. Suppose $X$ is an
$r$-normed quasi-Banach space.  Suppose $f\in
B_{\ell_{\infty}(E)}.$ Then there exist
$(c_j)_{j=1}^{\infty}$ and $(g_j)_{j=1}^{\infty}$ so that
$|c_j|\le M\delta^j$ and $f(\gamma)=\sum c_j\gamma(g_j)$ for
$\gamma\in E.$ Now if $\phi\in \Cal P_E(X)$ it follows that
$$ T_f\phi(h) =\sum_{j=1}^{\infty}c_j \phi(g_jh)$$ and so $$
\|T_f\phi\|_p \le
M(\sum_{j=1}^{\infty}\delta^{js})^{1/s}\|\phi\|_p$$ where
$s=\min(p,r).$ Thus $\|f\|_{\Cal M_p(E,X)}\le C$ where
$C=C(p,r,E)$ and so by Lemma 1 $E$ has property $\Cal
C_p(X).$

Now conversely suppose that $0<r<1$, $0<p\le\infty$ and that
$E$ has property $\Cal C_p(X)$ for every $r$-normable space
$X$.  It follows from consideration of
$\ell_{\infty}-$products that there exists a constant $C$ so
that for every $r$-normed space $X$ we have $\|f\|_{\Cal
M_p(E,X)}\le C\|f\|_{\infty}$ for $f\in\ell_{\infty}(E).$

Suppose $F$ is a finite subset of $E$.  We define an
$r$-norm $\|\ \|_A$ on $\ell_{\infty}(F)$ by setting $
\|f\|_A$ to be the infimum of $(\sum |c_j|^r)^{1/r}$ over
all $(c_j)_{j=1}^{\infty}$ and $(g_j)_{j=1}^{\infty}$ such
that $$ f(\gamma) =\sum_{j=1}^{\infty}c_j\gamma(g_j)$$ for
$\gamma\in F.$ Notice that $\|f_1f_2\|_A \le
\|f_1\|_A\|f_2\|_A$ for all $f_1,f_2\in A=\ell_{\infty}(F).$

For $\gamma\in F$ let $e_{\gamma}$ be defined by
$e_{\gamma}(\gamma)=1$ if $\gamma=\chi$ and $0$ otherwise.
Then for $f\in A,$ with $\|f\|_{\infty}\le 1,$ $$ (\int_G
\|\sum_{\gamma\in
F}f(\gamma)e_{\gamma}\gamma\|_A^pd\mu_G)^{1/p} \le C(\int_G
\|\sum_{\gamma\in F}e_{\gamma}\gamma\|_A^pd\mu_G)^{1/p}.$$
But for any $g\in G$ $\|\sum \gamma(g)e_{\gamma}\|_A \le 1.$
Define $H$ to be subset of $h\in G$ such that
$\|\sum_{\gamma\in F}f(\gamma)\gamma(h)e_{\gamma}\|_A \le
3^{1/p}C.$ Then $\mu_G(H)\ge 2/3$.  Thus there exist
$h_1,h_2\in H$ such that $h_1h_2=1$ (the identity in $G$).
Hence by the algebra property of the norm $$ \|f\|_A \le
3^{2/p}C^2 $$ and so if we fix an integer $C_0>3^{2/p}C^2$
we can find $c_j$ and $g_j$ so that $\sum |c_j|^r\le C_0^r$
and $$ f(\gamma)=\sum c_j\gamma(g_j)$$ for $\gamma\in F.$ We
can suppose $|c_j|$ is monotone decreasing and hence that
$|c_j|\le C_0j^{-1/r}.$ Choose $N_0$ so that
$C_0\sum_{j=N_0+1}^{\infty}j^{-1/r}\le 1/2.$ Thus $$
|f(\gamma) -\sum_{j=1}^{N_0}c_j\gamma(g_j)| \le 1/2$$ for
$\gamma\in F.$ Since each $|c_j|\le C_0$ this implies that
$B_{\ell_{\infty}(F)}\subset AP(F,N,1/2)$ where $N=C_0N_0.$

As this holds for every finite set $F$ it follows by an easy
compactness argument that $B_{\ell_{\infty}(E)}\subset
AP(E,N,1/2)$ and so by Theorem 2 $E$ is a set of
interpolation.\bull\enddemo

\proclaim{Theorem 4}Let $X$ be a natural quasi-Banach space
and suppose $0<p\le \infty.$ Then any Sidon set has property
$\Cal C_p(X).$\endproclaim

\demo{Proof} Suppose $E$ is a Sidon set.  Then there is a
constant $C_0$ so that if $f\in\ell_{\infty}(E)$ then there
exists $\nu\in C(G)^*$ such that $\hat
\mu(\gamma)=f(\gamma)$ for $\gamma\in E$ and $\|\mu\|\le
C_0\|f\|_{\infty}.$ We will show the existence of a constant
$C$ such that $\|f\|_{\Cal M_p(E,X)}\le C\|f\|_{\infty}.$ If
no such constant exists then we may find a sequence $E_n$ of
finite subsets of $E$ such that $\lim C_n=\infty$ where
$C_n$ is the least constant such that $\|f\|_{\Cal
M_p(E_n,X)}\le C_n\|f\|_{\infty}$ for all $f\in
\ell_{\infty}(E_n).$

Now the spaces $\Cal M_p(E_n,X)$ are each isometric to a
subspace of $\ell_{\infty}(L_p(G,X))$ and hence so is
$Y=c_0(\Cal M_p(E_n,X))$.  In particular $Y$ is natural.
Notice that $Y$ has a finite-dimensional Schauder
decomposition.  We will calculate the Banach envelope $Y_c$
of $Y.$ Clearly $Y_c=c_0(Y_n)$ where $Y_n$ is the
finite-dimensional space $\Cal M_p(E_n,X)$ equipped with its
the envelope norm $\|f\|_c.$

Suppose $f\in \ell_{\infty}(E_n)$.  Then clearly
$\|f\|_{\infty}\le \|f\|_{\Cal M_p(E,X)}$ and so
$\|f\|_{\infty}\le \|f\|_c.$ Conversely if
$f\in\ell_{\infty}(E_n)$ there exists $\nu\in C(G)^*$ with
$\|\nu\|\le C_0\|f\|_{\infty}$ and such that $\int \gamma
\,d\nu = f(\gamma)$ for $\gamma\in E_n.$ In particular
$C_0^{-1}\|f\|_{\infty}^{-1}f$ is in the absolutely closed
closed convex hull of the set of functions $\{\tilde g:g\in
G\}$ where $\tilde g(\gamma)=\gamma(g)$ for $\gamma\in E_n.$
Since $\|\tilde g\|_{\Cal M_p(E,X)}=1$ for all $g\in G$ we
see that $\|f\|_{\infty}\le \|f\|_c\le C_0\|f\|_{\infty}.$

This implies that $Y_c$ is isomorphic to $c_0.$ Since $Y$
has a finite-dimensional Schauder decomposition and is
natural we can apply Theorem 3.4 of [6] to deduce that
$Y=Y_c$ is already locally convex.  Thus there is a constant
$C_0'$ independent of $n$ so that $\|f\|_{\Cal M_p(E,X)}\le
C_0'\|f\|_{\infty}$ whenever $f\in\ell_{\infty}(E_n).$ This
contradicts the choice of $E_n$ and proves the
theorem.\bull\enddemo

We now consider the case of A-convex quasi-Banach spaces.
For this notion we will introduce the concept of an analytic
Sidon set.  We say a subset $E$ of $\Gamma$ is an {\it
analytic Sidon set} if $E$ satisfies $\Cal C_{\infty}(X)$
for every A-convex quasi-Banach space $X$.

\proclaim{Proposition 5}Suppose $0<p<\infty$.  Then $E$ is
an analytic Sidon set if and only if $E$ satisfies $\Cal
C_p(X)$ for every A-convex quasi-Banach space
$X.$\endproclaim

\demo{Proof}Suppose first $E$ is an analytic Sidon set, and
that $X$ is an A-convex quasi-Banach space (for which we
assume the quasinorm is plurisubharmonic).  Then $L_p(G,X)$
also has a plurisubharmonic quasinorm and so $E$ satisfies
(1) for $X$ replaced by $L_p(G,X)$ and $p$ replaced by
$\infty$ with constant $M$.  Now suppose $x_1,\ldots,x_n\in
X$ and $\gamma_1,\ldots,\gamma_n \in E.$ Define
$y_1,\ldots,y_n\in L_p(G,X)$ by $y_k(g)=\gamma_k(g)x_k.$
Then $$ \max_{g\in
G}\|\sum_{k=1}^ny_k\gamma_k(g)\|_{L_p(G,X)}=\|\sum_{k=1}^nx_k\gamma_k\|_p
$$ and a similar statement holds for the characters
$\epsilon_k$ on the Cantor group.  It follows quickly that
$E$ satisfies (1) for $p$ and $X$ with constant $M.$

For the converse direction suppose $E$ satisfies $\Cal
C_p(X)$ for every A-convex space $X$.  Suppose $X$ has a
plurisubharmonic quasinorm.  We show that $\Cal
M_{\infty}(E,X)=\ell_{\infty}(E).$ In fact $\Cal
M_{\infty}(F,X)$ can be isometrically embedded in
$\ell_{\infty}(X)$ for every finite subset $F$ of $E.$ Thus
(1) holds for $X$ replaced by $\Cal M_{\infty}(F,X)$ for
some constant $M$, independent of $F.$ Denoting by
$e_{\gamma}$ the canonical basis vectors in
$\ell_{\infty}(E)$ we see that if
$F=\{\gamma_1,\ldots,\gamma_n\}\subset E$ then $$
(\int_{\Delta}\|\sum_{k=1}^n\epsilon_k(t)e_{\gamma_k}\|^p_{\Cal
M_{\infty}(F,X)}d\mu_{\Delta}(t))^{1/p} \le M \max_{g\in
G}\|\sum_{k=1}^n \gamma_k(g)e_{\gamma_k}\|_{\Cal
M_{\infty}(F,X)}=M.$$

Thus the set $K$ of $t\in\Delta$ such that
$\|\sum_{k=1}^n\epsilon_k(t)e_{\gamma_k}\| \le 3^{1/p}M$ has
measure at least 2/3.  Arguing that $K.K=\Delta$ we obtain
that $$ \|\sum_{k=1}^n\epsilon_k(t)e_{\gamma_k}\|_{\Cal
M_{\infty}(F,X)}\le 3^{2/p}M^2$$ for every $t\in\Delta.$ It
follows quite simply that there is a constant $C$ so that
for every real valued $f\in\ell_{\infty}(F)$ we have
$\|f\|_{\Cal M_{\infty}(E,X)}\le C\|f\|_{\infty}.$ In fact
this is proved by writing each such $f$ with
$\|f\|_{\infty}=1$ in the form
$f(\gamma_k)=\sum_{j=1}^{\infty}2^{-j}\epsilon_k(t_j)$ for a
suitable sequence $t_j\in\Delta.$ A similar estimate for
complex $f$ follows by estimating real and imaginary parts.
Finally we conclude that since these estimates are
independent of $F$ that $\ell_{\infty}(E)=\Cal
M_{\infty}(E,X).$ \bull\enddemo

Of course any set of interpolation is an analytic Sidon set
and any analytic Sidon set is a Sidon set.  The next theorem
will show that not every analytic Sidon set is a set of
interpolation.  If we take $G=\bold T$ and $\Gamma=\bold Z$,
we recall that a Hadamard gap sequence is a sequence
$(\lambda_k)_{k=1}^{\infty}$ of positive integers such that
for some $q>1$ we have $\lambda_{k+1}/\lambda_k\ge q$ for
$k\ge 1.$ It is shown in [10] and [14] that a Hadamard gap
sequence is a set of interpolation.  However the union of
two such sequences may fail to be a set of interpolation;
for example $(3^n)_{n=1}^{\infty}\cup(3^n+n)_{n=1}^{\infty}$
is not a set of interpolation, since the closures of $(3^n)$
and $(3^n+n)$ in $b\bold Z$ are not disjoint.

\proclaim{Theorem 6}Let $G=\bold T$ so that $\Gamma=\bold
Z.$ Suppose $ E\subset\bold N$ is a finite union of Hadamard
gap sequences.  Then $E$ is an analytic Sidon
set.\endproclaim

\demo{Proof}Suppose $E=(\lambda_k)_{k=1}^{\infty}$ where
$(\lambda_k)$ is increasing.  We start with the observation
that $E$ is the union of $m$ Hadamard sequences if and only
there exists $q>1$ so that $\lambda_{m+k}\ge q^m\lambda_k$
for every $k\ge 1.$

We will prove the theorem by induction on $m.$ Note first
that if $m=1$ then $E$ is a Hadamard sequence and hence [14]
a set of interpolation.  Thus by Theorem 2 above, $E$ is an
analytic Sidon set.

Suppose now that $E$ is the union of $m$ Hadamard sequences
and that the theorem is proved for all unions of $l$
Hadamard sequences where $l<m.$ We assume that
$E=(\lambda_k)$ and that there exists $q>1$ such that
$\lambda_{k+m}\ge q^m\lambda_k$ for $k\ge 1.$ We first
decompose $E$ into at most $m$ Hadamard sequences.  To do
this let us define $E_1=\{\lambda_1\}\cup\{\lambda_k:k\ge
2,\ \lambda_k\ge q\lambda_{k-1}\}.$ We will write
$E_1=(\tau_k)_{k\ge 1}$ where $\tau_k$ is increasing.  Of
course $E_1$ is a Hadamard sequence.

For each $k$ let $D_k=E\cap [\tau_k,\tau_{k+1}).$ It is easy
to see that $|D_k|\le m$ for every $k$.  Further if $n_k\in
D_k$ then $n_{k+1}\ge \tau_{k+1}\ge qn_k$ so that $(n_k)$ is
a Hadamard sequence.  In particular $E_2=E\setminus E_1$ is
the union of at most $m-1$ Hadamard sequences and so $E_2$
is an analytic Sidon set by the inductive hypothesis.

Now suppose $w\in \bold T$.  We define
$f_w\in\ell_{\infty}(E)$ by $f_w(n)=w^{n-\tau_k}$ for $n\in
D_k.$ We will show that $f_w$ is uniformly continuous for
the Bohr topology on $\bold Z$; equivalently we show that
$f_w$ extends to a continuous function on the closure
$\tilde E$ of $E$ in the Bohr compactification $b\bold Z$ of
$\bold Z.$ Indeed, if this is not the case there exists
$\xi\in\tilde E$ and ultrafilters $\Cal U_0$ and $\Cal U_1$
on $E$ both converging to $\xi$ so that $\lim_{n\in\Cal
U_0}f_w(n) =\zeta_0$ and $\lim_{n\in\Cal U_1}f_w(n)=\zeta_1$
where $\zeta_1\neq \zeta_0.$ We will let
$\delta=\frac13|\zeta_1-\zeta_0|.$

We can partition $E$ into $m$ sets $A_1,\ldots,A_m$ so that
$|A_j\cap D_k|\le 1$ for each $k.$ Clearly $\Cal U_0$ and
$\Cal U_1$ each contain exactly one of these sets.  Let us
suppose $A_{j_0}\in\Cal U_0$ and $A_{j_1}\in\Cal U_1.$

Next define two ultrafilters $\Cal V_0$ and $\Cal V_1$ on
$\bold N.$ $\Cal V_0=\{V:\cup_{k\in V}D_k\in \Cal U_0\}$ and
$\Cal V_1=\{V:\cup_{k\in V}D_k\in\Cal U_1\}.$ We argue that
$\Cal V_0$ and $\Cal V_1$ coincide.  If not we can pick
$V\in\Cal V_0\setminus\Cal V_1$.  Consider the set
$A=(A_{j_0}\cap(\cup_{k\in V}D_k))\cup
(A_{j_1}\cap(\cup_{k\notin V}D_k)).$ Then $A$ is a Hadamard
sequence and hence a set of interpolation.  Thus for the
Bohr topology the sets $A_{j_0}\cap (\cup_{k\in V}D_k))$ and
$A_{j_1}\cap(\cup_{k\notin V}D_k))$ have disjoint closures.
This is contradiction since of course $\xi$ must be in the
closure of each.  Thus $\Cal V_0=\Cal V_1.$

Since both $\Cal U_0$ and $\Cal U_1$ converge to the same
limit for the Bohr topology we can find sets $H_0\in\Cal
U_0$ and $H_1\in\Cal U_1$ so that if $n_0\in H_0,\ n_1\in
H_1$ then $|w^{n_1}-w^{n_0}|<\delta$ and further
$|f_w(n_0)-\zeta_0|< \delta$ and
$|f_w(n_1)-\zeta_1|<\delta.$

Let $V_0=\{k\in\bold N:D_k\cap H_0\neq\emptyset\}$ and
$V_1=\{k\in\bold N:D_k\cap H_1\neq \emptyset\}.$ Then
$V_0\in \Cal V_0$ and $V_1\in\Cal V_1.$ Thus $V=V_0\cap
V_1\in \Cal V_0=\Cal V_1.$ If $k\in V$ there exists $n_0\in
D_k\cap H_0$ and $n_1\in D_k\cap H_1.$ Then $$ \align
3\delta &= |\zeta_1-\zeta_0|\\ &< |f_w(n_1)-f_w(n_0)| +
2\delta\\ &= |w^{n_1}-w^{n_0}| + 2\delta\\ &< 3\delta.
\endalign $$ This contradiction shows that each $f_w$ is
uniformly continuous for the Bohr topology.

Now suppose that $X$ is an $r$-normed A-convex quasi-Banach
space where the quasi-norm is plurisubharmonic.  Since both
$E_1$ and $E_2$ are analytic Sidon sets we can introduce a
constant $C$ so that if $f\in \ell_{\infty}(E_j)$ where
$j=1,2$ then $\|f\|_{\Cal M_{\infty}(E_j,X)}\le
C\|f\|_{\infty}.$ Pick a constant $0<\delta<1$ so that
$3.4^{1/r}\delta<C.$

Let $K_l=\{w\in\bold T:f_w\in AP(E,l,\delta)\}.$ It is easy
to see that each $K_l$ is closed and since each $f_w$ is
uniformly continuous by the Bohr topology it follows from
the Stone-Weierstrass theorem that $\cup K_l=\bold T.$ If we
pick $l_0$ so that $\mu_{\bold T}(K_{l_0})>1/2$ then
$K_{l_0}K_{l_0}=\bold T$ and hence since the map $w\to f_w$
is multiplicative $f_w\in AP(E,l_0^2,3\delta)$ for every
$w\in\bold T.$

Let $F$ be an arbitrary finite subset of $E.$ Then there is
a least constant $\beta$ so that $\|f\|_{\Cal
M_{\infty}(F,X)}\le \beta\|f\|_{\infty}.$ The proof is
completed by establishing a uniform bound on $\beta.$

For $w\in \bold T$ we can find $c_j$ with $|c_j|\le 1$ and
$\zeta_j\in \bold T$ for $1\le j\le l_0^2$ such that $$
|f_w(n) -\sum_{j=1}^{l_0^2}c_j\zeta_j^n|\le 3\delta$$ for
$n\in E.$ If $\tilde\zeta_j$ is define by
$\tilde\zeta_j(n)=\zeta_j^n$ then of course
$\|\tilde\zeta_j\|_{\Cal M_{\infty}(E,X)}=1.$ Restricting to
$F$ we see that $$ \|f_w\|_{\Cal M_{\infty}(F,X)}^r \le
l_0^2 + \beta^r(3\delta)^r.$$

Define $F:\bold C\to \Cal M_{\infty}(F,X)$ by $F(z)(n) =
z^{n-\tau_k}$ if $n\in D_k.$ Note that $F$ is a polynomial.
As in Theorem 5, $\Cal M_{\infty}(F,X)$ has a
plurisubharmonic norm.  Hence $$ \|F(0)\|^r \le
\max_{|w|=1}\|F(w)\|^r \le l_0^2 +(3\delta)^r\beta^r.$$
Thus, if $\chi_A$ is the characteristic function of $A,$ $$
\|\chi_{E_1\cap F}\|_{\Cal M_{\infty}(F,X)}^r\le l_0^2
+(3\delta)^r\beta^r.$$ It follows that $$ \|\chi_{E_2\cap
F}\|_{\Cal M_{\infty}(F,X)}^r \le l_0^2 + (3\delta)^r\beta^r
+1.$$

Now suppose $f\in \ell_{\infty}(F)$ and $\|f\|_{\infty}\le
1.$ Then $$ \|f\chi_{E_j\cap F}\|_{\Cal M_{\infty}(F,X)}\le
\|f\chi_{E_j\cap F}\|_{\Cal M_{\infty}(E_j\cap
F,X)}\|\chi_{E_j\cap F}\|_{\Cal M_{\infty}(F,X)}$$ for
$j=1,2.$ Thus $$ \|f\|_{\Cal M_{\infty}(F,X)}^r\le
C^r(1+2l_0^2+2(3\delta)^r\beta^r).$$ By maximizing over all
$f$ this implies $$ \beta^r \le C^r(1+2l_0^2 +
2(3\delta)^r\beta^r)$$ which gives an estimate $$ \beta^r
\le 2C^r(1+2l_0^2)$$ in view of the original choice of
$\delta.$ This estimate, which is independent of $F,$
implies that $E$ is an analytic Sidon set.\bull\enddemo

\demo{Remark}We know of no example of a Sidon set which is
not an analytic Sidon set.
\enddemo

\vskip1truecm

\subheading{References}

\item{1.}N.  Asmar and S.J.  Montgomery-Smith, On the
distribution of Sidon series, Arkiv Math. to appear.

\item{2.}D.  Grow, A class of $I_0$-sets, Colloq.  Math. 53
(1987) 111-124.

\item{3.}S.  Hartmann and C. Ryll-Nardzewski, Almost
periodic extensions of functions, Colloq.  Math. 12 (1964)
23-39.

\item{4.}S.  Hartmann and C. Ryll-Nardzewski, Almost
periodic extensions of functions, II, Colloq.  Math. 12
(1964) 79-86.

\item{5.}J.P.  Kahane, Ensembles de Ryll-Nardzewski et
ensembles de H. Helson, Colloq.  Math. 15 (1965) 87-92.

\item{6.}N.J.  Kalton, Banach envelopes of non-locally
convex spaces, Canad.  J. Math. 38 (1986) 65-86.

\item{7.}N.J.  Kalton, Plurisubharmonic functions on
quasi-Banach spaces, Studia Math. 84 (1986) 297-324.

\item{8.}J.F.  M\'ela, Sur les ensembles d'interpolation de
C. Ryll-Nardzewski et de S. Hartmann, Studia Math. 29 (1968)
167-193.

\item{9.}J.F.  M\'ela, Certains ensembles exceptionnels en
analyse de Fourier, Ann.  Inst.  Fourier (Grenoble) 18
(1968) 32-71.

\item{10.}J.  Mycielski, On a problem of interpolation by
periodic functions, Colloq.  Math. 8 (1961) 95-97.

\item{11.}A.  Pe\l czy\'nski, Commensurate sequences of
characters, Proc.  Amer.  Math.  Soc. 104 (1988) 525-531.

\item{12.}G.  Pisier, Les in\'egalit\'es de Kahane-Khintchin
d'apr\`es C. Borell, S\'eminaire sur la g\'eometrie des
\'espaces de Banach, Ecole Polytechnique, Palaiseau,
Expos\'e VII, 1977-78.

\item{13.}C.  Ryll-Nardzewski, Concerning almost periodic
extensions of functions, Colloq.  \newline Math. 12 (1964)
235-237.

\item{14.}E.  Strzelecki, Some theorems of interpolation by
periodic functions, Colloq.  Math. 12 (1964) 239-248.

\bye